\documentclass[a4paper,12pt]{amsart}
\usepackage{ifthen}
\usepackage{graphicx}
\usepackage{mathrsfs}
\numberwithin{equation}{section}
\nonstopmode
\newtheorem{thm}{Theorem}[section]

\newtheorem{lem}[thm]{Lemma}

\theoremstyle{definition}

\newenvironment{pf}[1][]{%
 \vskip 3mm
 \noindent
 \ifthenelse{\equal{#1}{}}%
  {{\slshape Proof. }}%
  {{\slshape #1.} }%
 }%
{\qed\bigskip}

\newcounter{alphabet}
\newcounter{tmp}
\newenvironment{Thm}[1][]{\refstepcounter{alphabet}%
\bigskip%
\noindent%
{\bf Theorem \Alph{alphabet}}%
\ifthenelse{\equal{#1}{}}{}{ (#1)}%
{\bf .}
\itshape}{\vskip 8pt}

\newcommand{\A}{{\mathcal A}}

\newcommand{\C}{{\mathbb C}}
\newcommand{\D}{{\mathbb D}}

\newcommand{\Hol}{{\mathcal O}}

\newcommand{\R}{{\mathbb R}}

\newcommand{\es}{{\mathcal S}}
\newcommand{\T}{{\mathcal T}}
\newcommand{\U}{{\mathcal U}}

\newcommand{\Z}{{\mathbb Z}}

\renewcommand{\Im}{\,{\operatorname{Im}\,}}

\renewcommand{\Re}{{\operatorname{Re}\,}}

\renewcommand{\mod}{{\operatorname{mod}\,}}

\renewcommand{\arg}{\,{\operatorname{arg}\,}}

\newcounter{minutes}
\divide\time by 60
\newcounter{hours}
\multiply\time by 60
\addtocounter{minutes}{-\time}

\begin{document}
\bibliographystyle{amsplain}
\title{
Univalent functions with half-integral coefficients
}

\begin{center}
{\tiny \texttt{FILE:~\jobname .tex,
        printed: \number\year-\number\month-\number\day, 
        \thehours.\ifnum\theminutes<10{0}\fi\theminutes}
}
\end{center}

\author[N.~Hiranuma]{Naoki Hiranuma}
\address{Sumitomo Life Insurance Company \\
7-18-24 Tsukiji, Chuo-ku, Tokyo 104-8430, Japan}
\email{hiranumanaoki@gmail.com}
\author[T. Sugawa]{Toshiyuki Sugawa}
\address{Graduate School of Information Sciences, Tohoku University \\
Aoba-ku, Sendai 980-8579, Japan}
\email{sugawa@math.is.tohoku.ac.jp}
\begin{abstract}
B.~Friedman found in his 1946 paper that the set of analytic univalent functions
on the unit disk in the complex plane with integral Taylor coefficients
consists of nine functions.
In the present paper, we prove that the similar set obtained by replacing
``integral" by ``half-integral" consists of another twelve functions
in addition to the nine.
We also observe geometric properties of the twelve functions.
\end{abstract}
\thanks{
The present research was supported in part by JSPS Grant-in-Aid for 
Scientific Research (B) 22340025.
}
\maketitle

\section{Introduction}
Let $\A$ denote the set of analytic functions $f$ on the unit disk
$\D=\{z\in\C: |z|<1\}$ normalized so that $f(0)=0, f'(0)=1.$
The set $\es$ of univalent functions in $\A$ has been a central object
to study in Geometric Function Theory since Bieberbach \cite{Bieb16} gave the
conjecture that $|a_n|\le n$ for $n=2,3,\dots$ for a function
$f(z)=z+a_2z^2+a_3z^3+\dots$ in $\es$ and that equality holds precisely when
$f$ is the Koebe function
$$
K(z)=\frac{z}{(1-z)^2}=\sum_{n=1}^\infty nz^n
=z+2z^2+3z^3+\cdots
$$
or its rotation $e^{-i\theta}K(e^{i\theta}z).$
The Bieberbach conjecture had been a driving force to develop
Geometric Function Theory for a long time and was finally solved
in the affirmative by de Branges in 1985.
We state it for later reference.

\begin{Thm}[de Branges]
Let $f(z)=z+a_2z^2+a_3z^3+\cdots$ be a function in $\es.$
Then $|a_n|\le n$ for each $n\ge2.$
If equality holds for some $n\ge2,$ then $f$ is a rotation of the Koebe function.
\end{Thm}

We remark that the assertion was verified earlier for first several $n$'s.
See, for example, \cite{Gong:BC} for a history and a proof of the Bieberbach
conjecture.
Meanwhile, Friedman \cite{Fried46} proved the following interesting theorem.

\begin{Thm}[Friedman]
Suppose that all the coefficients $a_n$ of a function $f$ in $\es$ are 
$($rational$)$ integers.
Then $f$ is one of the following nine functions:
$$
z,\hspace{1em}\frac{z}{1\pm z},\hspace{1em} \frac{z}{1\pm z^2},\hspace{1em}
\frac{z}{(1\pm z)^2},\hspace{1em}\frac{z}{1\pm z+z^2}.
$$
\end{Thm}

Note that $z/(1\pm z+z^2)=z(1\mp z)/(1\mp z^3).$
We observe that the above nine functions have the form
$z/P(z)$ for a polynomial $P(z)$ of degree at most $2.$
Indeed, the crucial point in the simplified proof by Linis \cite{Linis55}
is the fact that the coefficients of
the function $F(\zeta)=1/f(1/\zeta)=\zeta+b_0+b_1/\zeta+b_2/\zeta^2+\cdots$
for $f(z)=z+a_2z^2+a_3z^3+\cdots$ are given by
\begin{align*}
b_n=(-1)^{n+1}
\begin{vmatrix}
a_2 & 1   & 0  & \cdots & 0 \\
a_3 & a_2 & 1 & \cdots & 0 \\
a_4 & a_3 & a_2 & \cdots & 0\\
\hdotsfor{5} \\
a_{n+2} & a_{n+1} & a_n & \cdots & a_2
\end{vmatrix}.
\end{align*}
In particular, if $a_n$ are all integers, then $b_n$ are integers, too.
On the other hand, Gronwall's area theorem (see \cite{Pom:univ}) asserts that
\begin{equation}\label{eq:area}
\sum_{n=1}^\infty n|b_n|^2\le 1.
\end{equation}
Therefore, $b_n=0$ for $n>1$ and $b_1\in\{-1,0,1\}$ whenever $b_n$ are all integers.
In this way, we can conclude that $z/f(z)$ is a polynomial of degree
at most 2.
This idea can be used for a more general situation.
For extensions of Theorem B to integers in an imaginary quadratic 
number field, see Shah \cite{Shah51},
Townes \cite{Townes54}, Linis \cite{Linis55}, Bernardi \cite{Ber56} and
Royster \cite{Royster56}.
Moreover, Jenkins \cite{Jenkins87} determined all those functions $f\in\es$
for which the coefficients of $1/f(1/\zeta)$ are either rational half-integers
or half-integers in an imaginary quadratic number field.
Here and hereafter, a half-integer will mean the half of an integer.
Note therefore that an integer is a half-integer in our context.

It may be a natural question to ask what we can say if we replace
``integers" by ``half-integers" in the assumption of Theorem B.
In this case, however, we would only deduce that $2^{n+1}b_n$ is an integer
for each $n\ge1$ merely from the above observation.
Indeed, when $f(z)=z-z^2/2$ (see \S\ref{ss:1}), 
we have $b_n=2^{-n-1}$ for $n\ge0.$
Nevertheless, we have a finiteness result even in a more general situation.
For a subset $E$ of $\C,$ let $\A(E)$ denote the set of functions
$f(z)=z+a_2z^2+a_3z^3+\cdots$ in $\A$ such that $a_n\in E$ for all $n\ge2.$
Set $\es(E)=\es\cap\A(E).$
Denote by $\D(a,r)$ the disk $|z-a|<r$ in the complex plane $\C.$
If $E\cap\D(a,r_0)=\{a\}$ for every $a\in E$ and for some constant $r_0>0$
which is independent of the point $a,$
then we will say that $E$ is {\it uniformly discrete} (with bound $r_0$).

\begin{thm}\label{thm:finite}
Suppose that $E\subset\C$ is uniformly discrete.
Then $\es(E)$ consists of finitely many functions.
\end{thm}

For instance, $\es(\frac 1N\Z)$ is a finite set for every natural number $N,$
where $\frac1N\Z=\{n/N: n\in\Z\}.$
Note also that the ring $\Hol$ of integers in an imaginary quadratic number field
is uniformly discrete.
Therefore, we obtain finiteness also for the case $E=\frac 1c\Hol$
for a non-zero element $c$ in $\Hol.$
In these cases, we can say more.
Indeed, the following remarkable result is a special case of
Salem's theorem \cite[Theorem II]{Salem45}.

\begin{Thm}[Salem]
Let $\Hol$ be either the ring of rational integers or
the ring of integers in an imaginary quadratic number field
and let $c$ be a non-zero element in $\Hol.$
Then each function in $\es(\frac1c\Hol)$ is a rational function.
\end{Thm}

Hence, it is, in principle, only a matter of complexity 
to determine $\es(\frac1c\Hol).$
We may manage to do that for $\es(\frac12\Z).$

\begin{thm}\label{thm:main}
Suppose that all the coefficients $a_n$ of a function $f$ in $\es$ are 
half-integers.
Then $f$ is either one of the nine functions in Theorem B or else one of the
following twelve functions:
$$
z\pm\frac{z^2}2,\hspace{1em}
\frac{z(2\pm z)}{2(1\pm z)},\hspace{1em}
\frac{z(2\pm z^2)}{2(1\pm z^2)},\hspace{1em}
\frac{z(2\pm z)}{2(1-z^2)},\hspace{1em}
\frac{z(2\pm z)}{2(1\pm z)^2},\hspace{1em}
\frac{z(2\pm z+z^2)}{2(1\pm z+z^2)}.
$$
\end{thm}

Obradovi\'c and Ponnusamy \cite{OP01} pointed out that
the nine functions in Friedman's theorem are all starlike and belong to
the class $\U$ of functions $f\in\A$
satisfying the inequality $|z^2f'(z)/f(z)^2-1|<1$ on $|z|<1.$
We cannot, however, say the same for the additional twelve functions
in Theorem \ref{thm:main}.
Indeed, the function $f(z)=z(2+z+z^2)/2(1+z+z^2)$ is not even
close-to-convex as we will see in \S \ref{ss:4}.
Since $z^2f'(z)/f(z)^2-1=z^2(-1+2z+z^2)/(2+z+z^2)^2,$
we can see that $|z^2f'(z)/f(z)^2-1|^2$ takes the value $5+10\sqrt2/3$ for
the choice $z=(-1+\sqrt7 i)/\sqrt 8\in\partial\D.$
Hence, this function does not belong to $\U.$

Note here that, throughout the present paper, the fraction
$a/b\cdot c$ will mean $a/(bc)$ to reduce the use of parentheses.

We briefly describe the organization of the present paper.
In Section 2, we prepare necessary tools for the proof of main results
as well as a proof of Theorem \ref{thm:finite}.
In Section 3, we observe geometric properties of the twelve functions
in Theorem \ref{thm:main} as part of a proof of it.
Section 4 is devoted to the proof of Theorem \ref{thm:main}.
We also provide a collection of formulae which are useful in the proof.

\medskip
\noindent
\textbf{Acknowledgements.}
The authors would like to thank S.~Ponnusamy for bringing their attention
to a paper \cite{Jenkins87} by Jenkins.

\section{Necessary conditions for univalence}

For a function $f(z)=z+a_2z^2+a_3z^3+\cdots$ in $\A,$ we expand
$1/f$ in the Laurent series
$$
\frac 1{f(z)}=\frac1z+\sum_{n=0}^\infty b_nz^n
$$
on $0<|z|<\delta$ for sufficiently small $\delta>0.$
We denote by $\T$ the set of functions $f\in\A$ for which the inequality
\eqref{eq:area} holds.
Gronwall's area theorem means that $\es\subset\T.$
We set $\T(E)=\T\cap\A(E)$ for a subset $E$ of $\C.$

Following the idea of Friedman \cite{Fried46}, we now show the uniqueness lemma.

\begin{lem}\label{lem:unique}
Let $E$ be a uniformly discrete subset of $\C$ with bound $r_0$
and let $f(z)=z+a_2z^2+a_3z^3+\cdots$ and $g(z)=z+A_2z^2+A_3z^3+\cdots$
be functions in the class $\T(E).$
We write $1/f(z)=1/z+b_0+b_1z+\cdots.$
Suppose that $a_n=A_n$ for $n=2,\dots, N$ and that
\begin{equation}\label{eq:N}
2\sqrt{\frac{1-\sum_{n=1}^{N-2}n|b_n|^2}{N-1}} <r_0.
\end{equation}
Then $f=g.$
\end{lem}

\begin{pf}
We prove that $a_n=A_n$ for all $n$ by induction.
Assume that $a_n=A_n$ for $n=1,2,\dots, m$ with $m\ge N.$
By assumption, we have $f(z)=g(z)+cz^{m+1}+O(z^{m+2})$ as $z\to0,$
where $c=a_{m+1}-A_{m+1}.$
Hence, $f(z)/g(z)=1+cz^m+O(z^{m+1}),$ which leads to the expansion
\begin{align*}
\frac1{g(z)}&=\frac1{f(z)}\cdot\big(1+cz^m+O(z^{m+1})\big) \\
&=\frac1z+b_0+b_1z+\cdots+b_{m-2}z^{m-2}+(b_{m-1}+c)z^{m-1}+O(z^m).
\end{align*}
By noting that $g\in \T,$
we apply \eqref{eq:area} and \eqref{eq:N} to obtain
$$
(m-1)|b_{m-1}+c|^2
\le 1-\sum_{n=1}^{m-2}n|b_n|^2
\le 1-\sum_{n=1}^{N-2}n|b_n|^2
<\frac{(N-1)r_0^2}4
\le\frac{(m-1)r_0^2}4.
$$
Hence, $|b_{m-1}+c|<r_0/2.$
Similarly, the assumption that $f\in\T$ leads to the inequality
$|b_{m-1}|<r_0/2.$
The triangle inequality now yields $|a_{m+1}-A_{m+1}|=|c|<r_0.$
Since $E$ is uniformly discrete with bound $r_0,$
we obtain $a_{m+1}=A_{m+1}.$
By induction, we conclude that $a_n=A_n$ for all $n;$
namely, $f=g.$
\end{pf}

\begin{pf}[Proof of Theorem \ref{thm:finite}]
We are now ready to prove Theorem \ref{thm:finite}.
Suppose that $E$ is uniformly discrete with bound $r_0$
and let $N$ be a natural number so large
that $1<(N-1)r_0^2/4.$
We note that the condition \eqref{eq:N} is fulfilled whatever
$b_n$'s are.
Since $|a_n|\le n$ holds for all $n$ by the de Branges theorem
for $f(z)=z+a_2z^2+\cdots$ in $\es,$ we have only finitely
many choices of $a_2, a_3, \dots, a_N$ as the coefficients
of functions in $\es(E).$
Once $a_2, a_3, \dots, a_N$ are specified, by Lemma \ref{lem:unique},
there is at most one candidate for such a function $f\in\es(E).$
The proof is now complete.
\end{pf}

When $E=\frac12\Z,$ we can take $1/2$ as the bound $r_0.$
Therefore, the above proof tells us that
it is enough to examine all possible values of $a_2,\dots, a_{17}.$
By virtue of the de Branges theorem,
except for rotations of the Koebe function,
possible values of $a_n$ are $0, \pm 1/2, \pm 1, \dots, \pm (2n-1)/2.$
Therefore, without any additional constraint,
the number of these possibilities would be $7\cdot 11\cdots (4\cdot17-1)
\approx 8.14\times 10^{23}.$
To exclude non-univalent cases, we need more effective criteria for
univalence.

For a function $f\in\A,$ we expand the analytic function $\log(f(z)-f(w))/(z-w)$
in the polydisk $|z|<r, |w|<r$ for small enough $r>0$ in the form
\begin{align*}
\log\frac{f(z)-f(w)}{z-w}
&=-\sum_{j,k=0}^\infty c_{j,k}z^jw^k \\
&=-\sum_{j,k=1}^\infty c_{j,k}z^jw^k
+\log\frac{f(z)}{z}+\log\frac{f(w)}{w}.
\end{align*}

The coefficients $c_{j,k}$ are called Grunsky coefficients of $f.$
Grunsky's inequality was strengthened by Pommerenke \cite{Pom:univ}
as follows:
If $f\in\A$ is univalent on $|z|<1$ then
$$
\sum_{m=1}^n m\left|\sum_{k=1}^n c_{m,k}t_k\right|^2
\le\sum_{m=1}^\infty m\left|\sum_{k=1}^n c_{m,k}t_k\right|^2
\le\sum_{m=1}^n\frac{|t_m|^2}m
$$
for arbitrary $n\ge1$ and $t_1,\dots,t_n\in\C.$
This implies that the Hermitian matrix
$G_f(n)=(\gamma_{j,k}^{(n)})$ of order $n$ is positive semi-definite, where
$$
\gamma_{j,k}^{(n)}=\frac{\delta_{j,k}}j-\sum_{m=1}^nmc_{m,j}\overline{c_{m,k}}
$$
and $\delta_{j,k}$ means Kronecker's delta.
We will call $G_f(n)$ the Grunsky matrix of order $n$ for $f.$
We remark that $G_f(n)$ can be expressed in terms of $a_2,\dots,a_{2n+1}$
(see Appendix).

Since the above inequalities imply $|c_{j,k}|\le1/\sqrt{jk}\le 1,$
these are sufficient conditions for univalence, as well (see \cite{Pom:univ}).
We summarize these observations in the following form.

\begin{lem}\label{lem:gr}
A function $f\in\A$ is univalent on $\D$ if and only if
its Grunsky matrix $G_f(n)$ of order $n$
is positive semi-definite for every $n\ge1.$
\end{lem}

Prawitz's inequality, which is an extension of Gronwall's inequality,
is also useful as a univalence criterion.
See \cite{Milin:univ} for details.

\begin{lem}[Prawitz's inequality]\label{lem:prawitz}
Let $f\in \es$ and $[z/f(z)]^{\alpha}=1-\sum_{n=1}^{\infty}\sigma_nz^n$.
Then 
$$
\sum_{n=1}^{\infty}(n-\alpha )|\sigma_n|^2 \le\alpha
$$
for every $\alpha>0.$
\end{lem}

It is elementary, but not easy by hand, to compute the Grunsky matrices
for a given function.
However, by using a suitable computer software, we can check positivity
of $G_f(n)$ rigorously for a specific $f$ and a small enough $n.$
We collect useful formulae to compute these coefficients in Appendix.

\section{Properties of the twelve functions}

As part of the proof of Theorem \ref{thm:main},
we check univalence of the twelve functions in this section.
These functions may be in a special position within the class $\es.$
We will see geometric properties of these functions as well.
Since each pair can be interchanged by a suitable rotation,
it is enough to consider one function of each pair in Theorem \ref{thm:main}.

We recall here special classes of univalent functions.
See \cite{Pom:univ} as a fundamental reference.
A function $f\in\A$ is called {\it starlike} if $f$ maps $\D$
univalently onto a starlike domain with respect to the origin.
It is well known that $f\in\A$ is starlike if and only if
$\Re[zf'(z)/f(z)]>0$ on $|z|<1.$
For instance, the Koebe function is starlike.
A function $f\in\A$ is called {\it close-to-convex} if
$\Re[e^{i\theta}zf'(z)/g(z)]>0$ on $|z|<1$ for some $\theta\in\R$ and
a starlike function $g\in\es.$
Note that a starlike function is close-to-convex.
The Noshiro-Warschawski theorem implies that a close-to-convex function
is univalent.
Therefore, it is enough to show that $f$ is close-to-convex in order to check
univalence of $f.$
A more intrinsic characterization of close-to-convex functions was
given by Kaplan \cite{Kaplan52}.
For a locally univalent function $f\in\A$ we define
\begin{equation}\label{eq:F}
F_r(\theta)=\arg\left[\frac{\partial}{\partial\theta}f(re^{i\theta})\right]
=\arg f'(re^{i\theta})+\frac\pi2+\theta
\end{equation}
so that $F_r$ is continuous on $\R$ and satisfies the relation
$F_r(\theta+2\pi)=F_r(\theta)+2\pi.$
Then such an $f$ is close-to-convex if and only if
\begin{equation}\label{eq:kaplan}
F_r(\theta_2)-F_r(\theta_1)=\int_{\theta_1}^{\theta_2}
\Re\left[ 1+\frac{re^{i\theta}f''(re^{i\theta})}{f'(re^{i\theta})}\right]
d\theta>-\pi
\end{equation}
whenever $0<r<1$ and $\theta_1<\theta_2.$

\subsection{The function $f_1(z)=z+z^2/2$}\label{ss:1}
It is easy to check that $|zf_1'(z)-f_1(z)|\le|zf_1'(z)+f_1(z)|$
on $|z|<1,$ which is equivalent to $\Re[zf_1'(z)/f_1(z)]>0$ on $|z|<1.$
Hence, $f_1$ is starlike.
It is well known that $f_1$ maps $\D$ univalently
onto the inside of a cardioid.

\subsection{The function $f_2(z)=z(2-z)/2(1-z)$}
We first note that
$$
f_2(z)=\frac{z(1-z/2)}{1-z}
=z+\sum_{n=2}^\infty\frac{z^n}2
=z+\frac{z^2}2+\frac{z^3}2+\frac{z^4}2+\cdots.
$$
Since $p_2(z)=zf_2'(z)/f_2(z)=(2-2z+z^2)/(1-z)(2-z),$ a simple computation gives
$$
p_2(e^{i\theta})=\frac{i}{\theta}-\frac12+O(\theta)
$$
as $\theta\to0.$
By continuity, we can see that $\Re p(z)<0$ for a point $z\in\D$ close to $1.$
This means that $f$ is not starlike.
On the other hand, by taking the Koebe function $K(z),$ we have
$$
\Re\frac{zf_2'(z)}{K(z)}=\Re\left(1-z+\frac{z^2}2\right)>0,
$$
which implies that $f_2$ is close-to-convex and, therefore, univalent.

\subsection{The function $f_3(z)=z(2-z^2)/2(1-z^2)$}
This is expanded in the form
$$
f_3(z)=\frac{z(1-z^2/2)}{1-z^2}
=z+\sum_{n=2}^\infty\frac{z^{2n-1}}2
=z+\frac{z^3}2+\frac{z^5}2+\frac{z^7}2+\cdots.
$$
Since $p_3(z)=zf_3'(z)/f_3(z)=(2-z^2+z^4)/(1-z^2)(2-z^2)$ 
has the asymptotic behaviour
$p_3(e^{i\theta})=i/\theta-2+O(\theta)$ as $\theta\to0,$ 
the function $f_3$ is not starlike.
On the other hand, letting $g_3$ be the starlike function $z/(1-z^2),$
we have
$$
\frac{zf_3'(z)}{g_3(z)}=\frac{2-z^2+z^4}{2(1-z^2)}
=\frac12\left(\frac{1+z^2}{1-z^2}+1-z^2\right),
$$
which obviously has positive real part on $\D.$
Therefore, $f_3$ is close-to-convex.

\subsection{The function $f_4(z)=z(2+z)/2(1-z^2)$}
We can expand as follows:
$$
f_4(z)=\frac{z(1+z/2)}{1-z^2}
=\sum_{n=1}^\infty\left(z^{2n-1}+\frac{z^{2n}}{2}\right)
=z+\frac{z^2}2+z^3+\frac{z^4}2+z^5+\frac{z^6}2+\cdots.
$$
If we take the starlike function $g_4(z)=z/(1-z^2),$ then
$$
\frac{zf_4'(z)}{g_4(z)}=\frac{1+z+z^2}{1-z^2}
=\frac12\left(\frac{1+z}{1-z}+\frac{1+z^2}{1-z^2}\right),
$$
which has positive real part.
Therefore, $f_4$ is close-to-convex.
Since $4f_4(z)+1=(1+4z+z^2)/(1-z^2),$ we observe that
$$
f_4(e^{i\theta})=-\frac14-i\frac{2+\cos\theta}{4\sin\theta}.
$$
Therefore, $f_4$ maps $\D$ onto the complex plane slit along
the two half-lines $-1/4+iy,~|y|\ge\sqrt3/4.$
In particular, $f_4$ is not starlike.

\subsection{The function $f_5(z)=z(2-z)/2(1-z)^2$}
This function can be expressed by
$$
f_5(z)=\frac{z(1-z/2)}{(1-z)^2}
=\sum_{n=1}^\infty\frac{(n+1)z^n}2
=z+\frac{3z^2}2+2z^3+\frac{5z^4}2+3z^5+\frac{7z^6}2+\cdots.
$$
It is notable that the derivative has the simple form $f_5'(z)=1/(1-z)^3.$
In particular, $zf_5'(z)/K(z)=1/(1-z)$ has real part at least $1/2.$
Therefore, $f_5$ is close-to-convex.
It is easy to check that the boundary of the image $f_5(\D)$
is the parabola $x+2y^2+3/8=0.$
In particular, $f_5$ is not starlike, but it is
a concave function with opening angle $2\pi$
(see \cite{AW05} for its definition).

\subsection{The function $f_6(z)=z(2-z+z^2)/2(1-z+z^2)$}\label{ss:4}
This function and its rotation have the most complicated behaviour
among the twelve functions.
First $f_6$ can be expanded in the following form:
\begin{align*}
f_6(z)&=\frac{z(1+z/2+z^3/2)}{1+z^3}
=z+\frac{z^2}2+\sum_{n=1}^\infty \frac{(-1)^n}2\big(z^{3n+1}+z^{3n+2}\big) \\
&=z+\frac{z^2}2-\frac{z^4}2-\frac{z^5}2+\frac{z^7}2+\frac{z^8}2-\frac{z^{10}}2
-\frac{z^{11}}2+\frac{z^{13}}2+\frac{z^{14}}2-\frac{z^{16}}2-\cdots.
\end{align*}
We next observe that 
$f_6(e^{i\theta})=(2-e^{i\theta}+e^{2i\theta})/2(2\cos\theta-1).$
In particular, $2\Im f_6(e^{i\theta})=(\sin 2\theta-\sin\theta)/(2\cos\theta-1)
=\sin \theta.$
The denominator vanishes precisely when $\theta\equiv \pm\pi/3$
$(\mod 2\pi).$
We now show that $f_6$ is injective on the unit circle $\partial\D$
except for $e^{\pm\pi i/3}.$
Assume that $f_6(e^{is})=f_6(e^{it})\ne\infty$ for distinct points $e^{is}$
and $e^{it}$ in $\partial\D.$
By the symmetry of $f_6$ in $\R,$ we may assume that $0\le s\le\pi$
and $0\le t<2\pi.$
Then, by taking the imaginary part, we have $\sin s=\sin t,$ which
enforces $t=\pi-s.$
Therefore, letting $x=\cos t=-\cos s,$ we have
$\Re [f_6(e^{is})-f_6(e^{it})]=(1-x+2x^2)/(1-2x)-(1+x+2x^2)/(1+2x)
=2x(1+4x^2)/(1-4x^2)=0.$
This implies $x=0,$ which contradicts distinctness of the two points.
We have proved univalence of $f_6$ on the boundary of $\D$ except for
$e^{\pm\pi i/3}.$
We now apply (a slightly modified version of)
Darboux's theorem to ensure univalence of $f_6$ on $\D.$

We next show that $f_6$ is not close-to-convex.
Define $F_r(\theta)$ by \eqref{eq:F} for $f=f_6$ so that $F_r(0)=0.$
Let $F_1(\theta)$ be its limit as $r\to1-.$
Note that $\Re f_6(e^{i\theta})\to \pm\infty$ as $\theta\to\pi/3\mp$
whereas $\Im f_6(e^{i\theta})\to \sin(\pi/3)/2=\sqrt3/4$ as $\theta\to\pi/3.$
Therefore, $F_1(\theta)$ has a jump of $+2\pi$ at $\theta=\pi/3.$
By the symmetry in $\R,$ $F_1(\theta)$ has a jump of $+2\pi$ at $5\pi/3.$
Since $\frac{d}{d\theta}f_6(e^{i\theta})
=ie^{i\theta}(2-4\cos\theta+3\cos2\theta-i\sin2\theta)/2(2\cos\theta-1)^2,$
we have
$$
F_1(\theta)=\begin{cases}
q(\theta)+\theta+\pi/2 &\quad (0\le \theta<\pi/3), \\
q(\theta)+\theta+5\pi/2 &\quad (\pi/3<\theta<5\pi/3), \\
q(\theta)+\theta+9\pi/2 &\quad (5\pi/3<\theta<2\pi).
\end{cases}
$$
Here, $q(\theta)=\arg(2-4\cos\theta+3\cos2\theta-i\sin2\theta)$ is
the continuous branch determined by $q(0)=0.$
We note that $q(2\pi)=-4\pi.$
Observe that 
$F_1(0)=\pi/2, F_1(\pi/3-)=0, F_1(\pi/3+)=2\pi, F_1(\pi)=3\pi/2,
F_1(5\pi/3-)=\pi, F_1(5\pi/3+)=3\pi,$ and $F_1(2\pi)=5\pi/2.$
We now see that $2\Im f_6(e^{i\theta})=\sin\theta$ is increasing at
$\theta=\pi/3,$ which implies that $F_1(\pi/3+\delta)>2\pi$ for
small enough $\delta>0.$
By the symmetry, we also have $F_1(5\pi/3-\delta)<\pi$ for the same $\delta.$
Therefore, $F_1(5\pi/3-\delta)-F_1(\pi/3+\delta)<-\pi,$ which violates
condition \eqref{eq:kaplan}.
We have now proved that $f_6$ is not close-to-convex.

\subsection{Some pictures}
We present the images of $\D$ under these mappings, which are
generated by Mathematica 8.0.

\begin{figure}[htbp]
\begin{minipage}{0.3\hsize}
\begin{center}
\includegraphics[height=0.2\textheight]{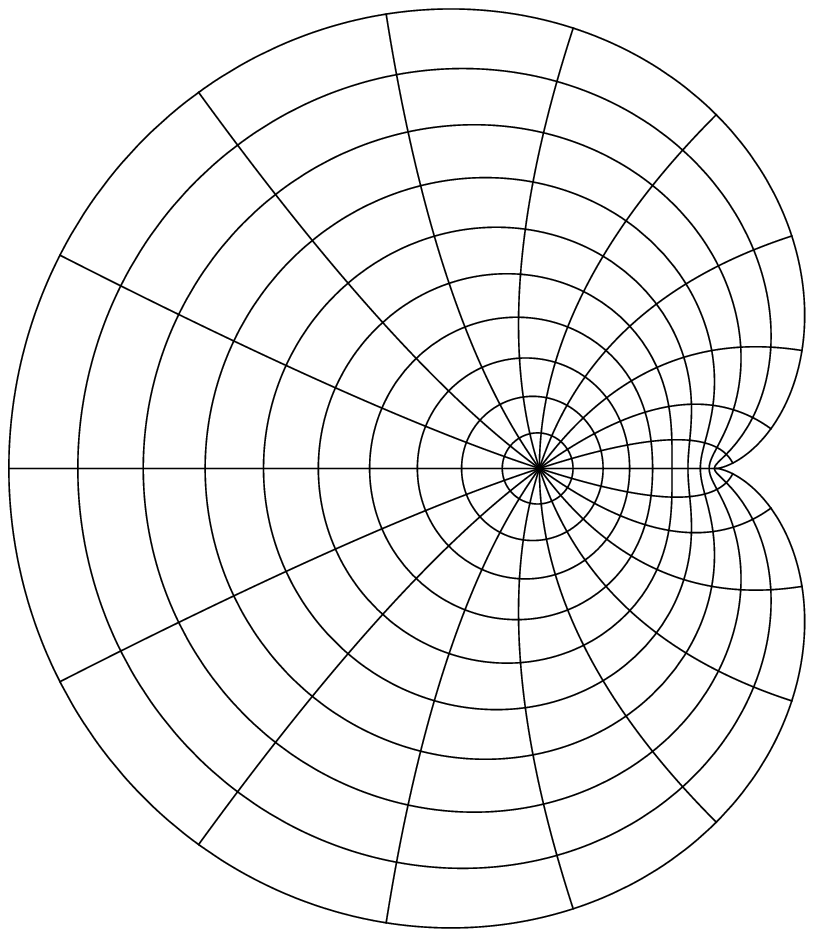}
\end{center}
\caption{$f_1(z)$}
\end{minipage}
\begin{minipage}{0.3\hsize}
\begin{center}
\includegraphics[height=0.2\textheight]{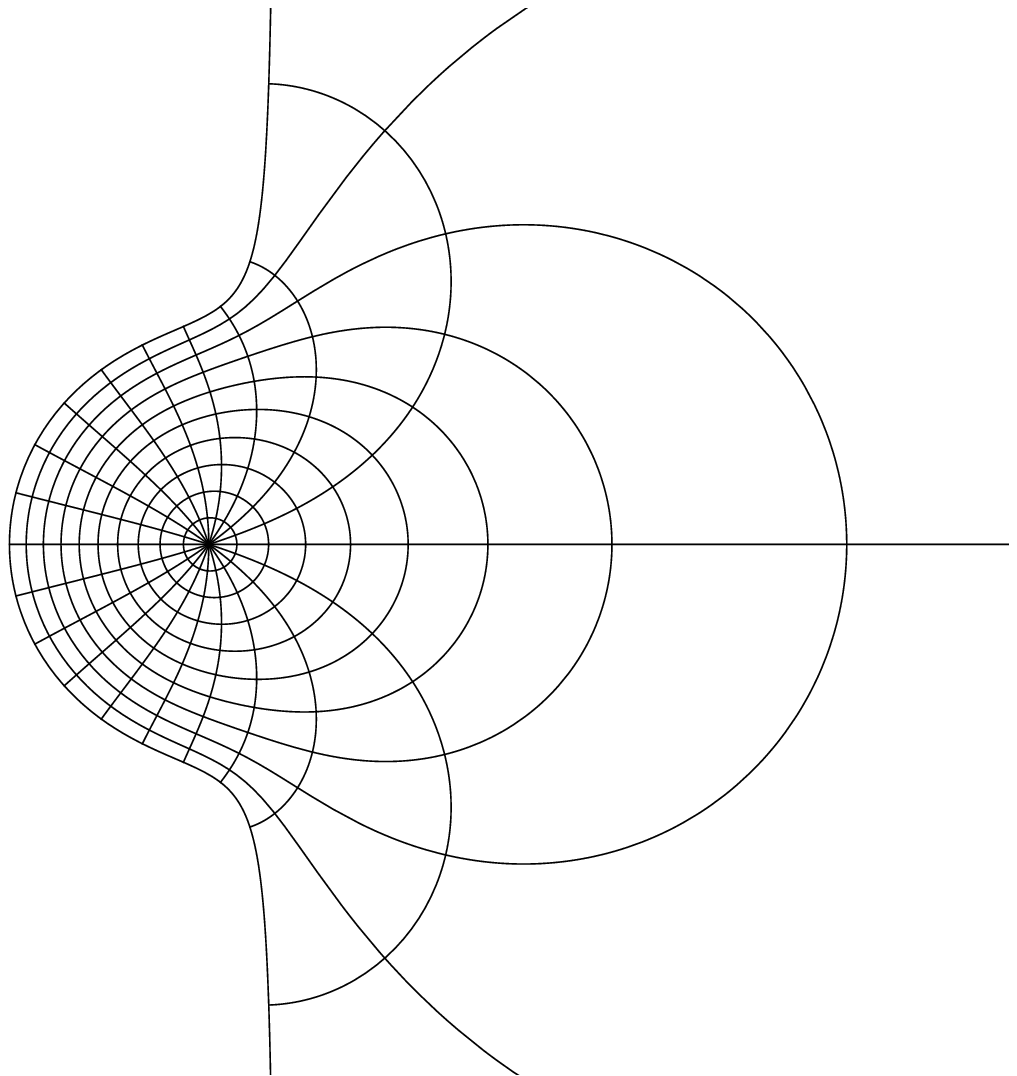}
\end{center}
\caption{$f_2(z)$}
\end{minipage}
\begin{minipage}{0.3\hsize}
\begin{center}
\includegraphics[height=0.2\textheight]{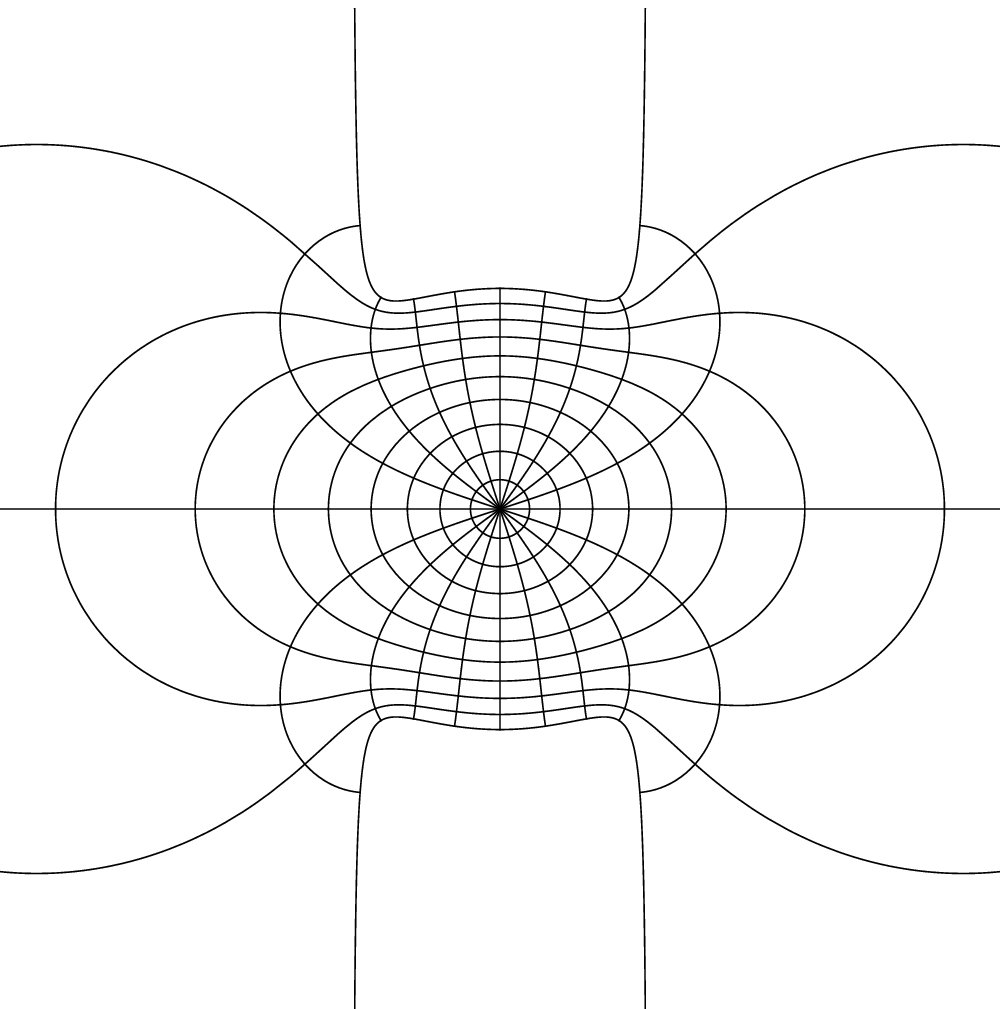}
\end{center}
\caption{$f_3(z)$}
\end{minipage}
\end{figure}

\begin{figure}[htbp]
\begin{minipage}{0.3\hsize}
\begin{center}
\includegraphics[height=0.2\textheight]{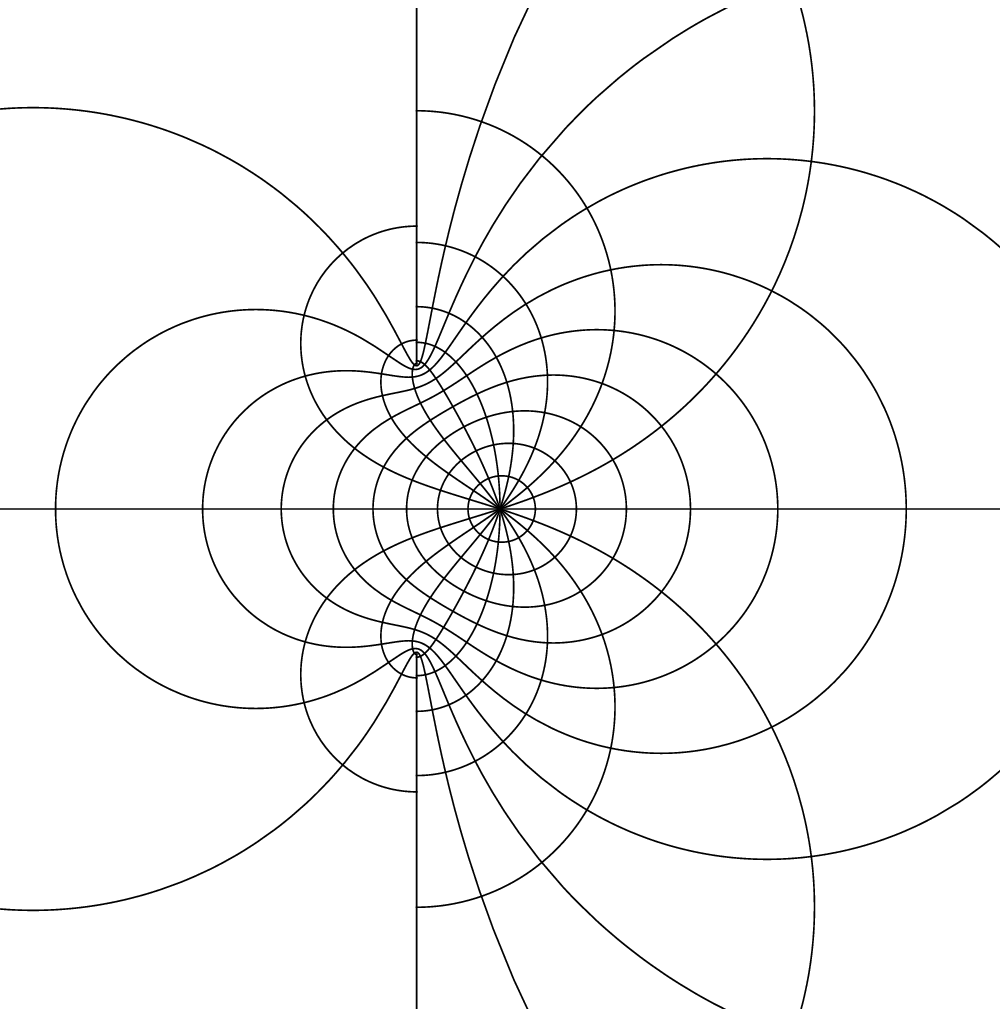}
\end{center}
\caption{$f_4(z)$}
\end{minipage}
\begin{minipage}{0.3\hsize}
\begin{center}
\includegraphics[height=0.2\textheight]{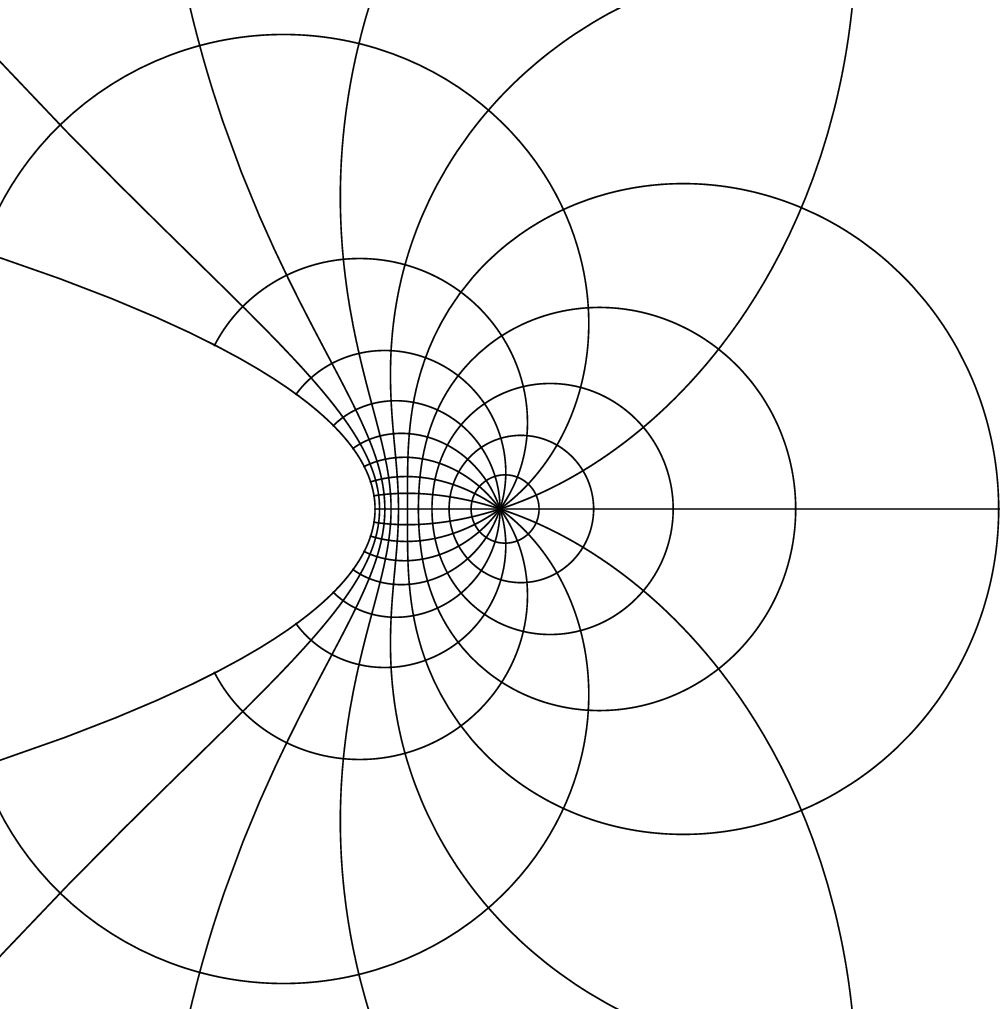}
\end{center}
\caption{$f_5(z)$}
\end{minipage}
\begin{minipage}{0.3\hsize}
\begin{center}
\includegraphics[height=0.2\textheight]{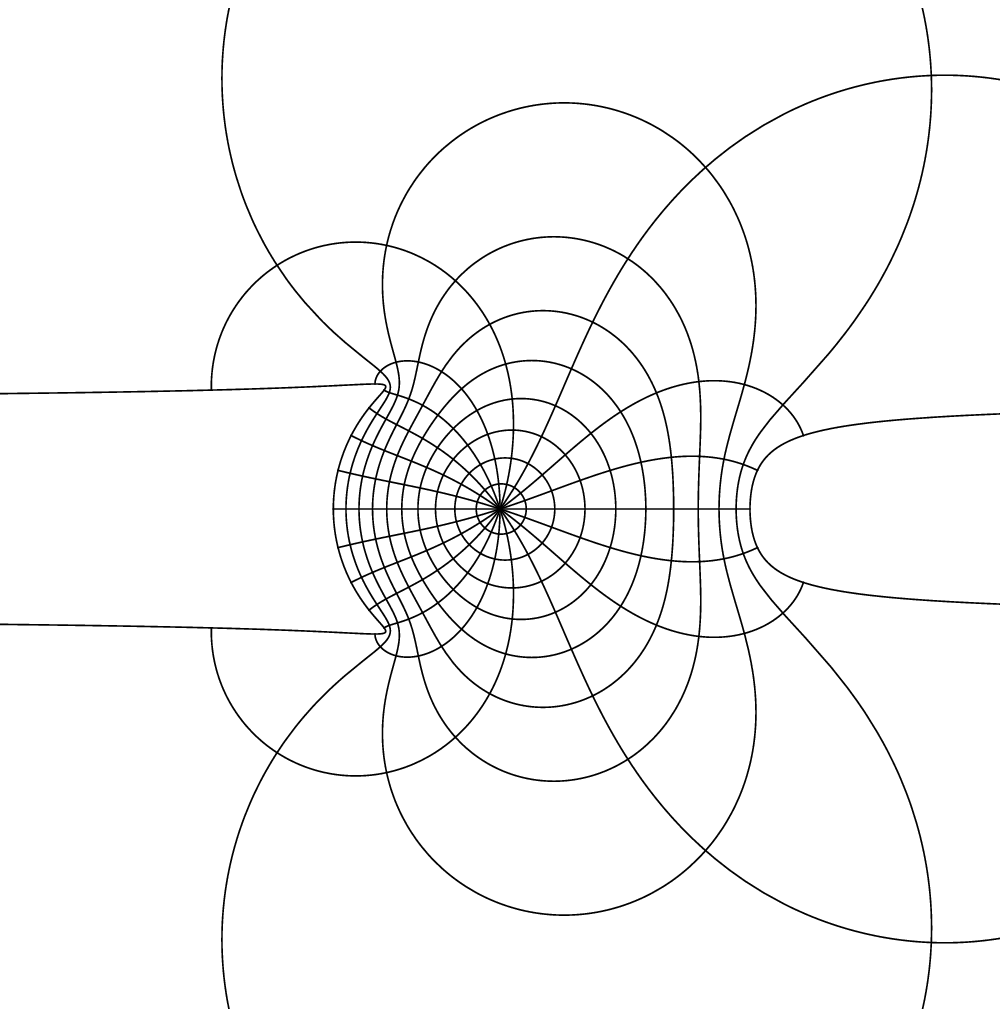}
\end{center}
\caption{$f_6(z)$}
\end{minipage}
\end{figure}

\section{Proof of Theorem \ref{thm:main}}

Let $f(z)=z+a_2z^2+a_3z^3+\cdots$ be a function in $\es(\frac12\Z)$
and fix it throughout the present section.
Then $2a_n$ is an integer for each $n.$
By taking the rotation $-f(-z)$ if necessary, we can assume that $a_2\ge0.$
Then $0\le a_2\le 2,$ and $a_2=2$ only when $f$ is the Koebe function by Theorem A.
Therefore, it suffices to consider the cases $a_2=0, 1/2, 1$ and $3/2.$
We will always assume that $f$ is not a rotation of the Koebe function
in the following so that $|a_n|<n$ holds for every $n\ge2$ by Theorem A.
As before, we write $1/f(z)=1/z+b_0+b_1z+b_2z^2+\cdots.$
For expressions of $b_n$ in terms of $a_n,$ see Appendix.
By \eqref{eq:ab},
if we specify $a_2,\dots,a_N,$ then $b_0,\dots, b_{N-2}$ are determined
the inequality
$$
|b_{N-1}|\le \sqrt{\frac{1-\sum_{n=1}^{N-2}n|b_n|^2}{N-1}}
$$
is obtained.
The inequality of this type will appear frequently in the sequel
without reference to the area theorem.
We remark that if the right-hand side is less than $1/4,$
Lemma \ref{lem:unique} guarantees that a function $f\in\es(\tfrac12\Z)$
with $f(z)=z+a_2z^2+\cdots+a_Nz^N+\cdots$ is uniquely determined.

\subsection{Case when $a_2=3/2$}
We start with the case when $a_2=3/2.$
Then we have $|b_1|=|a_3-9/4|\le 1,$
which is equivalent to $5/4\le a_3\le 13/4.$
Therefore, we have only the possibilities that $a_3=3/2, 2, 5/2.$
We will show that $a_3$ must be $2$ in this case.

Suppose that $a_3=3/2.$
Then $b_1=3/4$ and $b_2=-a_4+9/8.$
Since $|b_2|\le \sqrt{7/16\cdot 2}<0.47,$ $a_4=1, 3/2.$
Assume first that $a_4=1.$ 
Then $b_2=1/8, b_3=-a_5+3/16.$
Since $|b_3|\le \sqrt{13/32\cdot3}<0.37,$ $a_5=0,1/2.$
When $a_5=0,$ we have $\det G_f(2)=-215/8192<0.$
When $a_5=1/2,$ we have $b_3=-5/16$ and $b_4=-a_6+9/32.$
But there is no $a_6\in\frac12\Z$ such that $|b_4|\le\sqrt{29/256\cdot4}
<0.17.$
At any event, this case is discarded.

Suppose next that $a_3=5/2.$
Then $b_1=-1/4$ and $b_2=-a_4+33/8.$
We now have $|b_2|\le\sqrt{15/16\cdot 2}<0.69,$ which forces
$a_4=7/2,~b_2=5/8$ and $b_3=-a_5+79/16.$
Since $|b_3|\le\sqrt{5/32\cdot 3}<0.23,$ we have $4.7<a_5,$
which is not allowed.

Therefore, $a_3=2.$
Then $b_1=1/4$ and $b_2=-a_4+21/8.$
Since $|b_2|\le \sqrt{15/16\cdot 2}<0.69,$ we have the possibilities
$a_4=2,5/2,3.$
We will show that $a_4=5/2.$
If $a_4=2,$ we have $b_2=5/8$ and $b_3=-a_5+25/16.$
The condition $|b_3|\le\sqrt{5/32\cdot3}<0.23$ forces $a_5=3/2.$
Then $\det G_f(2)=-495/8192<0,$ which is a contradiction.
We now conclude that $a_4\ne2.$

If $a_4=3,$ we have $b_2=-3/8$ and $b_3=-a_5+73/16.$
Since $|b_3|\le\sqrt{21/32\cdot 3}<0.47,$ we have only the possibility
$a_5=9/2.$
Then $b_3=1/16$ and $b_4=-a_6+213/32.$
Since the condition $|b_4|\le\sqrt{165/256\cdot 4}<0.41$ implies
$a_6>6.26,$ this case does not occur.

We have proved that $a_4=5/2.$
Then $b_2=1/8$ and $b_3=-a_5+49/16.$
Since $|b_3|\le\sqrt{29/32\cdot3}<0.55,$ we have $a_5=3, 7/2.$
If $a_5$ were $7/2,$ then $b_3=-7/16$ and $b_4=-a_6+161/32.$
Since $|b_4|\le\sqrt{85/256\cdot4}<0.29,$ we obtain $a_6=5.$
It leads further to $b_4=1/32$ and $b_5=-a_7+457/64.$
Since $|b_5|\le\sqrt{21/64\cdot 5}<0.26,$ we have $a_7=7,$
which is excluded by assumption.

Hence, we have shown that $a_5=3$ in this case.
Then $b_3=1/16$ and $b_4=-a_6+113/32.$
Since $|b_4|\le\sqrt{229/256\cdot 4}<0.48,$ we have $a_6=7/2, 4.$
If $a_6$ was $4,$ then $b_4=-15/32$ and $b_5=-a_7+353/64.$
Since $|b_5|\le\sqrt{1/64\cdot 5}<0.06,$ we have $a_7=11/2.$
Then $b_5=1/64,~b_6=-a_8+977/128$ and
$|b_6|\le\sqrt{59/4096\cdot 6}<0.05$ implies $7.58<a_8<7.69,$
which is impossible.
Hence, $a_6=7/2.$
Then $b_4=1/32$ and $b_5=-a_7+4+1/64.$
Since $|b_5|\le\sqrt{57/64\cdot 5}<0.43,$ we have $a_7=4.$
Then $b_5=1/64$ and $b_6=-a_8+9/2+1/128.$
Since $|b_6|\le\sqrt{3643/4096\cdot 6}<0.39,$ we have $a_8=9/2.$
We can continue this process to obtain $a_n=(n+1)/2$ up to $n=17.$
Then Lemma \ref{lem:unique} implies that $f=f_5$ in this case.

\subsection{Case when $a_2=1$}
We recall that $b_1=-a_3+1$ in this case.
Since $|-a_3+1|=|b_1|\le1,$ we have the possibilities
$a_3=0,1/2,1,3/2,2.$
When $a_3=0,$ we have $b_1=1,$ which forces $b_n$ to be $0$
for all $n>1.$
Therefore, we have $1/f(z)=1/z-1+z;$ namely, $f(z)=z/(1-z+z^2),$
which appears in Theorem B.

First, we show that $a_3\ne2.$
If $a_3=2,$ we have $b_1=-1.$
In this case, we have similarly $f(z)=z/(1-z-z^2)=z+z^2+2z^3+3z^4+5z^5
+8z^6+\cdots.$
This function, however, is not in $\es$ because it violates 
the assertion of Theorem A.
(Its coefficients are known as Fibonacci numbers.)
The remaining three cases will be discussed in the following subsections.

\subsubsection{Case when $a_3=1/2$}
We will show that this case does not occur.
If we assume $a_3=1/2,$ we have $b_1=1/2$ and $b_2=-a_4.$
The condition $|b_2|\le\sqrt{3/4\cdot 2}<0.62$ now implies
$a_4=-1/2, 0, 1/2.$

Suppose first that $a_4=-1/2.$
Then $b_2=1/2$ and $b_3=-a_5-5/4.$
The condition $|b_3|\le\sqrt{1/4\cdot 3}<0.29$ enforces
$a_5=-3/2, -1.$
If $a_5=-3/2,$ then $\gamma_{2,2}^{(2)}=-1/32<0.$
Therefore, $a_5\ne-3/2.$
On the other hand, if $a_5=-1,$ then $b_3=-1/4$ and $b_4=-a_6-3/4.$
The condition $|b_4|\le\sqrt{1/16\cdot 4}=1/8$ means
$-7/8\le a_6\le -5/8,$ which is not allowed to hold.
In a similar way, we can exclude the case $a_4=1/2.$

Finally, we suppose that $a_4=0.$
Then $b_2=0$ and $b_3=-a_5-1/4.$
Since $|b_3|\le\sqrt{3/4\cdot 3}=1/2,$ we have $a_5=-1/2, 0.$
Suppose that $a_5=-1/2.$
Then $b_3=1/4$ and $b_4=-a_6-3/4.$
Since $|b_4|\le\sqrt{9/16\cdot 4}=3/8,$ we have $a_6=-1, -1/2.$
When $a_6=-1,$ we have $b_4=1/4, b_5=-a_7-9/8$ and $|b_5|\le\sqrt{5/16\cdot 5}
=1/4.$
Thus $a_7=-1$ and $\det G_f(3)=-11/256<0,$ which cannot happen.
When $a_6=-1/2,$ we have similarly $a_7=0$ and $\det G_f(3)=-11/256<0.$
Therefore, $a_5\ne-1/2.$
In the same way, we can show that $a_5\ne0.$
We have proved that $a_3\ne1/2$ if $a_2=1.$

\subsubsection{Case when $a_3=1$}
In this case, we have $b_1=0$ and $b_2=-a_4+1.$
The condition $|b_2|\le\sqrt{1/2}<0.71$ implies
$a_4=1/2, 1, 3/2.$
We will show that $a_4=1$ is the only admissible case.
Indeed, we first suppose that $a_4=1/2.$
Then $b_2=1/2$ and $b_3=-a_5-3/4.$
Since $|b_3|\le\sqrt{1/4\cdot3}<0.62,$ we have $a_5=-1,-1/2.$
If $a_5=-1,$ then $\gamma_{2,2}^{(2)}=-7/4<0.$
If $a_5=-1/2,$ then $\gamma_{2,2}^{(2)}=-1/4<0.$
At any choice, we are led to a contradiction.
Therefore, we conclude that $a_4\ne1/2.$

We next suppose that $a_4=3/2.$
Then $b_2=-1/2$ and $b_3=-a_5+2.$
Since $|b_3|\le\sqrt{1/2\cdot 3}<0.41,$ we have $a_5=2.$
Then $b_3=0$ and $b_4=-a_6+5/2.$
Since $|b_4|\le\sqrt{1/2\cdot 4}<0.36,$ we have $a_6=5/2.$
Then $b_4=0$ and $b_5=-a_7+13/4.$
Since $|b_5|\le\sqrt{1/2\cdot5}<0.32,$ we have $a_7=3, 7/2.$
When $a_7=3,$ we have $b_5=1/4$ and $b_6=-a_8+15/4.$
Since no element $a_8\in\frac12\Z$ satisfies 
the condition $|b_6|\le\sqrt{3/16\cdot6}<0.18,$ we see that $a_7\ne 3.$
We can similarly see that $a_7\ne7/2.$
At any event, the assumption $a_4=3/2$ yields a contradiction.
Thus we have seen that $a_4\ne3/2.$

Hence, we have shown that $a_4=1.$
Therefore, $b_2=0$ and $b_3=-a_5+1.$
Since $|b_3|\le \sqrt{1/3}<0.58,$ we have $a_5=1/2,1,3/2.$
Then $a_5=1$ is only the possible case.
Indeed, if $a_5=1/2,$ then $b_3=1/2$ and $b_4=-a_6.$
Since $|b_4|\le\sqrt{1/4\cdot4}=1/4,$ we have $a_6=0.$
Then $b_4=0$ and $b_5=-a_7-1/2.$
Since $|b_5|\le\sqrt{1/4\cdot5}<0.23,$ we have $a_7=-1/2.$
Then $b_5=0$ and $b_6=-a_8-1.$
Since $|b_6|\le\sqrt{1/4\cdot6}<0.21,$ we have $a_8=-1.$
Then $b_6=0$ and $b_7=-a_9-5/4.$
We now see that no element $a_9$ in $\frac12\Z$ does not satisfy
the condition $|b_7|\le\sqrt{1/4\cdot7}<0.19.$
Therefore, $a_5\ne1/2.$
Similarly, we can show that $a_5\ne3/2.$

Hence, we have shown that $a_5=1.$
Then $b_3=0$ and $b_4=-a_6+1.$
Since $|b_4|\le\sqrt{1/4}=1/2,$ we have $a_6=1/2,1,3/2.$
We will show that $a_6=1$ is the only possible case.
Indeed, if $a_6=1\pm1/2$ then $b_4=\mp1/2$ and $b_n=0$ for $n>4.$
Therefore, $1/f(z)=1/z-1\mp z^4/2$ and thus
$$
f(z)=\frac{z}{1-(z\pm z^5/2)}
=\sum_{n=0}^\infty z\left(z\pm \frac{z^5}2\right)^n.
$$
Therefore, $f\notin\A(\frac12\Z)$ (see $a_{11}$ for example).
Therefore, $a_6\ne\pm1/2.$
In this way, we obtain $a_6=1.$
Then $b_4=0$ and $b_5=-a_7+1.$
Since $|b_5|\le\sqrt{1/5}<0.45,$ we have $a_7=1.$
We can continue this process to obtain $a_n=1$ for $n\le 17.$

Here, we note that the function $z/(1-z)$ which appears in Theorem B
satisfies that $a_n=1$ for $1\le n\le 7.$
Lemma \ref{lem:unique} now implies that such a function is unique.
Hence, we have shown that $f(z)=z/(1-z)$ in this case.

\subsubsection{Case when $a_3=3/2$}
We will show that this case does not occur.
When $a_3=3/2,$ we have $b_1=-1/2$ and $b_2=-a_4+2.$
The condition $|b_2|\le\sqrt{3/4\cdot 2}<0.62$ implies
$a_4=3/2,2,5/2.$
Suppose first that $a_4=3/2.$
Then $b_2=1/2$ and $b_3=-a_5+7/4.$
Since $|b_3|\le\sqrt{1/4\cdot3}<0.29,$ we have $a_5=3/2, 2.$
If $a_5=3/2,$ then $\gamma_{2,2}^{(2)}=-1/32<0.$
If $a_5=2,$ then $\det G_f(2)=-11/128<0.$
Therefore, at any event, this case is not admitted.
The case when $a_4=5/2$ can be discarded in the same way.
We finally suppose that $a_4=2.$
Then $b_2=0$ and $b_3=-a_5+11/4.$
Since $|b_3|\le\sqrt{3/4\cdot3}=1/2,$ we have $a_5=5/2,3.$
If $a_5=5/2,$ then $b_3=1/4$ and $b_4=-a_6+13/4.$
Since $|b_4|\le\sqrt{9/16\cdot4}=3/8,$ we have $a_6=3,7/2.$
If $a_6=3$ in addition, then $b_4=1/4,~b_5=-a_7+29/8$
and $|b_5|\le\sqrt{5/16\cdot5}=1/4.$
In particular, $a_7=7/2$ and $\det G_f(3)=-1/128<0.$
Therefore, $a_6\ne3.$
In the same way, we have $a_6\ne7/2.$
Consequently, we have $a_5\ne5/2$ and thus $a_5=3.$
Then $b_3=-1/4$ and $b_4=-a_6+17/4.$
Since $|b_4|\le\sqrt{9/16\cdot4}=3/8,$ we have $a_6=4, 9/2.$
If $a_6=4,$ then $b_4=1/4$ and $b_5=-a_7+45/8.$
Since $|b_5|\le\sqrt{5/16\cdot5}=1/4,$ we have $a_7=11/2$
and $\det G_f(3)=-9/2048<0.$
The case when $a_6=9/2$ can also be discarded in the same way.
Therefore, $a_5\ne3$ and thus the possibility of $a_4=2$ has
been eliminated.

\subsection{Case when $a_2=1/2$}
When $a_2=1/2,$ $b_1=-a_3+1/4.$
Since $|b_1|\le1,$ we have $a_3=-1/2,0,1/2,1.$
We discuss these four cases separately.

\subsubsection{Case when $a_3=-1/2$}
We will show that this case is not allowed to occur.
Under the assumptions, we have $b_1=3/4$ and $b_2=-a_4-5/8.$
Since $|b_2|\le\sqrt{7/16\cdot2}<0.47,$ we have $a_4=-1,-1/2.$
If $a_4=-1,$ then $b_2=3/8$ and $b_3=-a_5-5/16.$
Since $|b_3|\le\sqrt{5/32\cdot3}<0.23,$ we have $a_5=-1/2.$
Then $\gamma_{2,2}^{(2)}=-41/512<0.$
This case is thus impossible.
If $a_4=-1/2,$ then $b_2=-1/8$ and $b_3=-a_5+3/16.$
Since $b_3|\le\sqrt{13/32\cdot3}<0.37,$ we have $a_5=0,1/2.$
When $a_5=0,$ we have $\det G_f(2)=-215/8192<0.$
When $a_5=1/2,$ we have $b_3=-5/16, b_4=-a_6+23/32$
and $|b_4|\sqrt{29/256\cdot4}<0.17.$
The last inequality implies $0.54<a_6<0.89,$ which is
satisfied by no $a_6$ in $\frac12\Z.$
Therefore, we have seen that this case does not occur.

\subsubsection{Case when $a_3=0$}
Then $b_1=1/4$ and $b_2=-a_4-1/8.$
Since $|b_2|\le\sqrt{15/16\cdot 2}<0.69,$ we have $a_4=-1/2,0,1/2.$
We will first show that the case $a_4=1/2$ does not occur.
Suppose, to the contrary, that $a_4=1/2.$
Then $b_2=-5/8$ and $b_3=-a_5+9/16.$
Since $|b_3|\le\sqrt{5/32\cdot3}<0.23,$ we have $a_5=1/2.$
We now have $\det G_f(2)=-495/8192<0,$ which is impossible.
Therefore, $a_4\ne1/2$ in this case.

We next suppose that $a_4=-1/2.$
Then $b_2=-5/8$ and $b_3=-a_5-7/16.$
Since $|b_3|\le\sqrt{21/32\cdot3}<0.47,$ we have $a_5=-1/2, 0.$
If $a_5=0,$ then $\det G_f(2)=-207/8192<0.$
Therefore, we must have $a_5=-1/2.$
Then $b_3=1/16$ and $b_4=-a_6-5/32.$
Since $|b_4|\le\sqrt{165/256\cdot4}<0.41,$ we have $a_6=-1/2, 0.$
If $a_6=-1/2,$ then $b_4=11/32$ and $b_5=-a_7-7/64.$
Since $|b_5|\le\sqrt{11/64\cdot5}<0.19,$ we have $a_7=0.$
Then $\gamma_{2,2}^{(3)}=-119/512<0.$
Therefore, we must have $a_6=0.$
In this case, $b_4=-5/32$ and $b_5=-a_7+25/64.$
Since $|b_5|\le\sqrt{35/64\cdot5}<0.34,$ we have $a_7=1/2.$
Then $b_5=-7/64$ and $b_6=-a_8+67/128.$
Since $|b_6|\le\sqrt{1995/4096\cdot6}<0.29,$ we have $a_8=1/2.$
Then $b_6=3/128$ and $b_7=-a_9+17/256.$
Since $|b_7|\le\sqrt{3963/8192\cdot6}<0.27,$ we have $a_9=0.$
Then $b_7=17/256$ and $b_8=-a_{10}-245/512.$
Since $|b_8|\le\sqrt{29681/65536\cdot8}<0.24,$ we have $a_{10}=-1/2.$
Here, we note that the function $f_6$ has the same coefficients so far.
Therefore, we now apply Lemma \ref{lem:unique} to conclude that
$f=f_6$ in this case.

Finally, we suppose that $a_4=0.$
Then $b_2=-1/8$ and $ b_3=-a_5+1/16.$
Since $|b_3|\le\sqrt{29/32\cdot3}<0.55,$ we have $a_5=0,1/2.$
If $a_5=1/2,$ we have $b_3=-7/16$ and $b_4=-a_6+15/32.$
Since $|b_4|\le\sqrt{85/256\cdot4}<0.29,$ we have $a_6=1/2.$
Then $b_4=-1/32$ and $b_5=-a_7+9/64.$
Since $|b_5|\le\sqrt{21/64\cdot5}<0.26,$ we have $a_7=0.$
Then $b_5=-9/64$ and $b_6=-a_8-17/128.$
Since $|b_6|\le\sqrt{939/4096\cdot6}<0.2,$ we have $a_8=0.$
Then $b_6=-17/128$ and $b_7=-a_9+89/256.$
Since $|b_7|\le\sqrt{1011/8192\cdot7}<0.14,$ we have
$0.2<a_9<0.49,$ which is impossible.
Therefore, the case when $a_5=1/2$ is discarded.
We thus have $a_5=0,$ which implies that $b_3=1/16$ and 
$b_4=-a_6-1/32.$
Since $|b_4|\le\sqrt{229/256\cdot4}<0.48,$ we have $a_6=-1/2, 0.$
If $a_6=-1/2,$ then $b_4=15/32,~b_5=-a_7-31/64$ and
$|b_5|\le\sqrt{1/64\cdot5}<0.06.$
The last inequality forces $a_7=-1/2$ and 
$\gamma_{2,2}^{(3)}=-55/512<0,$ which is not allowed.
Therefore, we must have $a_6=0.$
Then $b_4=-1/32$ and $b_5=-a_7+1/64.$
Since $|b_5|\le\sqrt{57/64\cdot5}<0.43,$ we have $a_7=0.$
We can continue this process to obtain $a_8=\cdots=a_{17}=0.$
Note here that the function $f_1(z)=z+z^2/2$ satisfies the above conditions.
We now apply Lemma \ref{lem:unique} to conclude that
$f=f_1$ in this case.

\subsubsection{Case when $a_3=1/2$}
Then $b_1=-1/4$ and $b_2=-a_4+3/8.$
Since $|b_2|\le\sqrt{15/16\cdot2}<0.69,$ we have $a_4=0,1/2,1.$
We will show that the only possible case is when $a_4=1/2.$

Indeed, we first suppose that $a_4=0.$
Then $b_2=3/8$ and $b_3=-a_5-1/16.$
Since $|b_3|\le\sqrt{21/32\cdot3}<0.47,$ we have $a_5=-1/2,0.$
If $a_5=-1/2,$ then $\gamma_{2,2}^{(2)}=-41/512<0.$
We thus must have $a_5=0.$
Then $b_3=-1/16$ and $b_4=-a_6-5/32.$
Since $|b_4|\le\sqrt{165/256\cdot4}<0.41,$ we have $a_6=-1/2,0.$
If $a_6=-1/2,$ then $b_4=11/32$ and $b_5=-a_7-25/64.$
Since $|b_5|\le\sqrt{11/64\cdot5}<0.19,$ we have $a_7=-1/2.$
Then $\det G_f(3)=-395595/2^{25}<0,$ which is impossible.
If $a_6=0,$ then $b_4=-5/32$ and $b_5=-a_7+7/64.$
Since $|b_5|\le\sqrt{35/64\cdot5}<0.34,$ we have $a_7=0.$
Then $\det G_f(3)=-955083/2^{25}<0,$ which is impossible, too.
Therefore, $a_4\ne0.$

We next suppose that $a_4=1.$
Then $b_2=-5/8$ and $b_3=-a_5+15/16.$
Since $|b_3|\le\sqrt{5/32\cdot3}<0.23,$ we have $a_5=1.$
Then $\det G_f(2)=-175/8192<0,$ which is impossible.
Therefore, $a_4\ne1.$

Hence, the remaining case is only when $a_4=1/2.$
In this case, $b_2=-1/8$ and $b_3=-a_5+7/16.$
Since $|b_3|\le\sqrt{29/32\cdot3}<0.55,$ we have $a_5=0, 1/2.$
We show that $a_5=1/2.$
If $a_5=0,$ then $b_3=7/16,~b_4=-a_6-1/32$ and 
$|b_4|\le\sqrt{85/256\cdot4}<0.29.$
Thus $a_6=0,~b_4=-1/32$ and $b_5=-a_7-9/64.$
Since $|b_5|\le\sqrt{21/64\cdot5}<0.26,$ we have $a_7=0.$
Then $\gamma_{3,3}^{(3)}=-31/1024<0,$ which is impossible.
Therefore, $a_5\ne0$ and thus $a_5=1/2.$
Then, $b_3=-1/16$ and $b_4=-a_6+15/32.$
Since $|b_4|\le\sqrt{229/256\cdot4}<0.48,$ we have $a_6=0,1/2.$
If $a_6=0,$ then $b_4=15/32$ and $b_5=-a_7-1/64.$
Since $|b_5|\le\sqrt{1/64\cdot5}<0.06,$ we have $a_7=0.$
Then $\gamma_{2,2}^{(3)}=-137/512<0,$ which is not allowed.
Therefore, we must have $a_6=1/2.$
Then $b_4=-1/32$ and $b_5=-a_7+31/64.$
Since $|b_5|\le\sqrt{57/64\cdot5}<0.43,$ we have $a_7=1/2.$
In the same way, we can show that $a_8=\cdots=a_{17}=1/2.$
Here, we note that the function $f_2$ has the same coefficients
as $f.$
Lemma \ref{lem:unique} now yields that $f=f_2$ in this case.

\subsubsection{Case when $a_3=1$}
Then $b_1=-3/4$ and $b_2=-a_4+7/8.$
Since $|b_2|\le\sqrt{7/16\cdot2}<0.47,$ we have $a_4=1/2, 1.$
We first show that $a_4\ne1.$
Suppose, to the contrary, that $a_4=1.$
Then $b_2=-1/8$ and $b_3=-a_5+21/16.$
Since $|b_3|\le\sqrt{13/32\cdot3}<0.37,$ we have $a_5=1,3/2.$
If $a_5=1,$ then $\gamma_{2,2}^{(2)}=-113/512<0.$
If $a_5=3/2,$ then $b_3=-3/16$ and $b_4=-a_6+55/32.$
Since $|b_4|\le\sqrt{77/256\cdot4}<0.28,$ we have $a_6=3/2.$
Then $b_4=7/32$ and $b_5=-a_7+133/64.$
Since $|b_5|\le\sqrt{7/64\cdot5}<0.15,$ we have $a_7=2.$
Then $\det G_f(3)=-523697/3\cdot2^{25}<0,$ which is impossible.
Therefore, both cases were discarded.

In this way, we have confirmed that $a_4\ne1$ and thus $a_4=1/2.$
Then $b_2=3/8$ and $b_3=-a_5+13/16.$
Since $|b_3|\le\sqrt{5/32\cdot3}<0.23,$ we have $a_5=1.$
We note that \eqref{eq:N} is satisfied with $N=4$
and that the univalent function $f_4$ has the same coefficients so far.
We now apply Lemma \ref{lem:unique} to conclude that $f=f_4$ in this case.

\subsection{Case when $a_2=0$}
Finally, we treat the case when $a_2=0.$
In this case, we have $b_1=-a_3.$
Since $|b_1|\le1,$ we have $a_3=0,\pm1/2, \pm1.$
When $a_3=1,$ then $b_n=0$ for $n>1$ by the area theorem.
Therefore, $1/f(z)=1/z-z;$ namely, $f(z)=z/(1-z^2),$ which appears in Theorem B.
When $a_3=-1,$ in the same way, we have $f(z)=z/(1+z^2).$
Therefore, we may restrict ourselves on the cases $a_3=-1/2, 0, 1/2.$
We shall consider each case in the following subsections.
Since the two cases when $a_3=\pm1/2$ can be treated similarly, we consider
these first.

\subsubsection{Case when $a_3=-1/2$}\label{sss:a3}
In this case, $b_1=1/2$ and $b_2=-a_4.$
Since $|b_2|\le\sqrt{3/4\cdot2}<0.62,$ we have $a_4=-1/2,0,1/2.$
We first show that $a_4\ne-1/2.$
Suppose, to the contrary, that $a_4=-1/2.$
Then $b_2=1/2$ and $b_3=-a_5+1/4.$
Since $|b_3|\le\sqrt{1/4\cdot3}<0.29,$ we have $a_5=0,1/2.$
If $a_5=0,$ then $\gamma_{2,2}^{(2)}=-1/32<0.$
If $a_5=1/2,$ then $b_3=-1/4$ and $b_4=-a_6+1/2.$
Since $|b_4|\le\sqrt{1/16\cdot4}=1/8,$ we have $a_6=1/2.$
Then $b_4=0$ and $b_5=-a_7-1/8.$
Since $|b_5|\le\sqrt{1/16\cdot5}<0.12,$ we have $-0.25<a_7<-0.005,$
which is impossible.
Therefore, we conclude that $a_4\ne-1/2.$

Next, we show that $a_4\ne1/2.$
To the contrary, suppose that $a_4=1/2.$
Then $b_2=-1/2$ and $b_3=-a_5+1/4.$
Since $|b_3|\le\sqrt{1/4\cdot3}<0.29,$ we have $a_5=0,1/2.$
If $a_5=0,$ then $\gamma_{2,2}^{(2)}=-1/32<0,$ which is impossible.
If $a_5=1/2,$ then $b_3=-1/4$ and $b_4=-a_6-1/2.$
Since $|b_4|\le\sqrt{1/16\cdot4}=1/8,$ we have $a_6=1/2.$
Then $b_4=0,~b_5=-a_7-1/8$ and $|b_5|\le\sqrt{1/16\cdot5}<0.12<1/8.$
There is no $a_7\in\frac12\Z$ in this case.
Hence, we have seen that $a_4\ne1/2.$

We now have only the possibility $a_4=0.$
Then, $b_2=0$ and $b_3=-a_5+1/4.$
Since $|b_3|\le\sqrt{3/4\cdot3}=1/2,$ we have $a_5=0,1/2.$
We show now that $a_5\ne0.$
If $a_5=0,$ then $b_3=1/4$ and $b_4=-a_6.$
Since $|b_4|\le\sqrt{9/16\cdot4}=3/8,$ we have $a_6=0.$
Then $b_4=0$ and $b_5=-a_7+1/8.$
Since $|b_5|\le\sqrt{9/16\cdot5}<0.34,$ we have $a_7=0.$
Then $\det G_f(3)=-49/2048<0,$ which is impossible.
Therefore, $a_5\ne0.$
We next assume that $a_5=1/2.$
Then $b_3=-1/4$ and $b_4=-a_6.$
Since $|b_4|\le\sqrt{9/16\cdot 4}=3/8,$ we have $a_6=0.$
Then $b_4=0$ and $b_5=-a_7-3/8.$
Since $|b_5|\le\sqrt{9/16\cdot5}<0.34,$ we have $a_7=-1/2.$
Then $b_5=1/8$ and $b_6=-a_8.$
Since $|b_6|\le\sqrt{31/64\cdot6}<0.29,$ we have $a_8=0.$
Then $b_6=0$ and $b_7=-a_9+7/16.$
Since $|b_7|\le\sqrt{31/64\cdot7}<0.26,$ we have $a_9=1/2.$
Then $b_7=-1/16$ and $b_8=-a_{10}.$
Since $|b_8|\le\sqrt{117/256\cdot8}<0.24,$ we have $a_{10}=0.$
Moreover, by Lemma \ref{lem:unique}, we have functions in $\es(\frac12\Z)$
with the above coefficients up to $n=10$ at most one.
We note here that the function $f_3(iz)/i$ has the above coefficients.
Therefore, we conclude that $f(z)=f_3(iz)/i$ in this case.

\subsubsection{Case when $a_3=1/2$}
In this case, $b_1=-1/2$ and $b_2=-a_4.$
Thus, as in the previous case, we have the possibilities
$a_4=-1/2,0,1/2.$
We first show that $a_4\ne-1/2.$
Indeed, if $a_4=-1/2,$ then $b_2=1/2$ and $b_3=-a_5+1/4.$
Since $|b_2|\le\sqrt{1/4\cdot2}<0.29,$ we have $a_5=0,1/2.$
When $a_5=0,$ we have $\gamma_{2,2}^{(2)}=-1/32<0.$
When $a_5=1/2,$ we have $\det G_f(2)=-11/128<0.$
At any event, the case when $a_4=-1/2$ is discarded.
In the same way, we can show $a_4\ne1/2.$

We have thus only the possibility that $a_4=0.$
In this case, as in \S\ref{sss:a3}, we conclude that $f=f_3.$

\subsubsection{Case when $a_3=0$}
In this case, we have $b_1=0,~b_2=-a_4,~b_3=-a_5,~b_4=-a_6.$
Since $|b_2|\le\sqrt{1/2}<0.71,$ we have $a_4=-1/2,0,1/2.$
We will show that $a_4\ne\pm1/2.$

To the contrary, we first suppose that $a_4=-1/2.$
Then $b_2=1/2$ and the condition $|b_3|\le\sqrt{1/2\cdot3}<0.41$
implies $a_5=0.$
Similarly, we further obtain $a_6=0$ and $b_5=-a_7+1/4.$
Since $|b_5|\le\sqrt{1/2\cdot5}<0.32,$ we have $a_7=0,1/2.$
If $a_7=0,$ then $\gamma_{3,3}^{(3)}=-5/12<0.$
Therefore, $a_7=1/2$ must hold.
Then $b_5=-1/4,~b_6=-a_8.$
Since $|b_6|\le\sqrt{3/16\cdot6}<0.18,$ we have $a_8=0$
and Lemma \ref{lem:unique} is applicable.
We now look at the function 
$$
g(z)=\frac{z(2+z^3)}{2(1+z^3)}=z-\frac{z^4}{2(1+z^3)}
=z-\frac{z^4}2+\frac{z^7}2-\frac{z^{10}}2+\cdots.
$$
Observe that $g(z)$ has the same coefficients as those of $f(z)$ up to $n=8$
and that
$$
\frac1{g(z)}=\frac1z+\frac{z^2}{2+z^3}
=\frac1z-\sum_{n=1}^\infty\frac{(-1)^n}{2^n}z^{3n-1}.
$$
Since $\sum_{n=1}^\infty (3n-1)2^{-2n}=1,$ we see that
$g\in\T(\frac12\Z).$
Therefore, Lemma \ref{lem:unique} implies that $f=g.$
However, $g$ is not univalent.
Indeed, we expand in the form
$$
\left[\frac z{g(z)}\right]^{2/3}=1+\sum_{n=1}^\infty\sigma_nz^n
=1+\frac{z^3}3-\frac{7z^6}{36}
+\frac{19z^9}{162}-\frac{143z^{12}}{1944}+\frac{281z^{15}}{5832}+O(z^{18}).
$$
Then we see that
$$
\frac23-\sum_{n=1}^{15}(n-\tfrac23)|\sigma_n|^2
=-\frac{353917}{2^6\cdot 3^{13}}<0.
$$
Namely, $g$ does not satisfy Prawitz's inequality (Lemma \ref{lem:prawitz})
with $\alpha=2/3.$
Therefore, this case is also discarded.
We have confirmed that $a_4\ne-1/2$ as long as $a_2=a_3=0.$
In the same way, we can show that $a_4\ne1/2.$

Hence, we have shown that $a_4=0.$
Then $1/f(z)=1/z-a_5z^3-a_6z^4-a_7z^5-a_8z^6+(a_5^2-a_9)z^7+O(z^8).$
Since $|b_3|\le\sqrt{1/3}<0.58,$ we have $a_5=-1/2,0,1/2.$
If $a_5=-1/2,$ we have $a_6=a_7=a_8=0$
and $|b_7|\le\sqrt{1/4\cdot7}<0.19,$ which implies $0.06<a_9<0.44.$
Therefore, we have $a_5\ne-1/2.$
Similarly, we have $a_5\ne1/2.$
Hence, we have shown that $a_5=0.$
Then $b_3=0$ and $|b_4|\le\sqrt{1/4}=1/2,$ which implies
$a_6=-1/2,0,1/2.$
If $a_6=-1/2,$ then $b_4=1/2$ and $b_n=0$ for $n>4.$
Therefore, $1/f(z)=1/z+z^4/2;$ namely,
$f(z)=z/(1+z^5/2).$
Then $f'(z)=4(1-2z^5)/(2+z^5),$ which has a zero in $|z|<1.$
This is impossible. Therefore, $a_6\ne-1/2.$
In the same way, we have $a_6\ne1/2.$

Hence, $a_6=0.$
Then $b_4=0$ and $|a_7|=|b_5|\le\sqrt{1/5}<1/2.$
Therefore, $a_7=0.$
We can continue this process to obtain $a_8=\dots=a_{17}=0.$
Lemma \ref{lem:unique} now implies that $f$ must be the
identity map; that is, $f(z)=z.$

The proof is now complete.

\section{Appendix}
In the present section, we collect several formulae which are 
useful in the proof of Theorem \ref{thm:main}.
Let $f(z)=z+a_2z^2+a_3z^3+\cdots$ be in $\A$ and
$1/f(z)=1/z+b_0+b_1z+b_2z^2+\cdots.$
We first note that the coefficients $b_n$ are computed
in terms of $a_n$'s recursively by the formula
\begin{equation}\label{eq:ab}
b_{n-1}=-a_{n+1}-\sum_{k=2}^n a_kb_{n-k},\quad n\ge1.
\end{equation}
In particular, we have
\begin{align*}
b_0&=-a_2, \\
b_1&=-a_3+a_2^2, \\
b_2&=-a_4+2a_2a_3-a_2^3, \\
b_3&=-a_5+2a_2a_4+a_3^2-3a_2^2a_3+a_2^4, \\
b_4&=-a_6+2a_2a_5+2a_3a_4-3a_2^2a_4-3a_2a_3^2+4a_2^3a_3-a_2^5,
\end{align*}
and so on.
The Grunsky coefficients $c_{j,k}$ of $f$ can be computed recursively by
$$
c_{j,k}=\sum_{l=1}^{k-1}\frac lk a_{k-l}c_{j+1,l}-\sum_{m=1}^j a_{m+1}c_{j-m,k}
-\frac{a_{j+k+1}}k
$$
for $j\ge0$ and $k\ge1$ (see \cite{KS12univ} for details).
Here, we set $a_1=1.$
It is easy to see that $c_{j,k}$ can be expressed as a polynomial
in $a_2,\dots,a_{j+k+1}.$
We also note that $c_{j,k}=c_{k,j}.$
For convenience, we write down the coefficients $c_{j,k}$ for $1\le j\le k\le3$
so that the reader can compute the Grunsky matrices of orders $2$ and $3:$
\begin{align*}
c_{1,1}&=-a_3+a_2^2, \\
c_{1,2}&=-a_4+2a_2a_3-a_2^3, \\
c_{1,3}&=-a_5+2a_2a_4+a_3^2-3a_2^2a_3+a_2^4, \\
c_{2,2}&=-a_5+2a_2a_4+\frac32a_3^2-4a_2^2a_3+\frac32a_2^4, \\
c_{2,3}&=-a_6+2a_2a_5+3a_3a_4-4a_2^2a_4-5a_2a_3^2+7a_2^3a_3-2a_2^5, \\
c_{3,3}&=-a_7+2a_2a_6+3a_3a_5-4a_2^2a_5+2a_4^2-12a_2a_3a_4 \\
&\qquad +8a_2^3a_4-\frac73a_3^3
+15a_2^2a_3^2-14a_2^4a_3+\frac{10}3a_2^6.
\end{align*}

\def\cprime{$'$} \def\cprime{$'$} \def\cprime{$'$}
\providecommand{\bysame}{\leavevmode\hbox to3em{\hrulefill}\thinspace}
\providecommand{\MR}{\relax\ifhmode\unskip\space\fi MR }
\providecommand{\MRhref}[2]{%
  \href{http://www.ams.org/mathscinet-getitem?mr=#1}{#2}
}
\providecommand{\href}[2]{#2}

\end{document}